\documentclass[12pt]{amsart}

\usepackage{latexsym}
\usepackage{amsmath}
\usepackage{amssymb}
\usepackage{amscd}
\usepackage{multicol}

\newtheorem{theorem}{\bf Theorem}[section]
\newtheorem{corollary}[theorem]{\bf Corollary}
\newtheorem{proposition}[theorem]{\bf Proposition}
\newtheorem{lemma}[theorem]{\bf Lemma}

\newtheorem{definition-theorem}[theorem]{\bf Theorem-Definition}

\def\bR{\mathbb{R}}
\def\bC{\mathbb{C}}
\def\bZ{\mathbb{Z}}

\def\D{{\mathcal D}}
\def\t{\mathfrak{t}}

\def\t{\frak{t}}
\def\p{\frak{p}}

\def\b{\frak{b}}
\def\u{\frak{u}}

\def\t{\frak{t}}

\def\b{\frak{b}}
\def\u{\frak{u}}

\def\r1{\sqrt{-1}}
\def\me{[[e^{t_i}]]}
\def\R{{\mathcal R}}
\def\S{{\mathcal S}}
\def\A{{\mathcal A}}

\def\hbar{h}

\def\bZ{\mathord{\Bbb Z}}
\def\bR{\mathord{\Bbb R}}

\def\bC{\mathord{\Bbb C}}

\def\t{\frak{t}}
\def\p1{\frac{\partial}{\partial t_1}}
\def\pi{\frac{\partial}{\partial t_i}}

\def\b{\frak{b}}
\def\u{\frak{u}}

\def\bZ{\mathord{\Bbb Z}}
\def\bR{\mathord{\Bbb R}}

\def\bC{\mathord{\Bbb C}}

\def\H{{\mathcal H}}

\def\H{{\mathcal H}}
\title{Relations in the quantum cohomology ring of $G/B$}

\date{\today}
\begin{document}
\begin{abstract}
The ideal of relations in the quantum cohomology of the flag
manifold $G/B$ has been determined by B. Kim in \cite{K}. We are
going to point out a limited number of properties that, if they
are satisfied by an $\bR[q_1,\ldots ,q_l]$-linear product $\circ$
on $H^*(G/B)\otimes \bR[q_1, \ldots, q_l]$, then the ring
$(H^*(G/B)\otimes \bR[q_1, \ldots ,q_l], \circ)$ is isomorphic to
Kim's ring.
\end{abstract}

\author[A.-L. Mare]{Augustin-Liviu  Mare}
\address{Department of Mathematics\\ University of Toronto
 \\Toronto, Ontario M5S 3G3, Canada}
 \email{amare@math.toronto.edu}

\maketitle

\section{Introduction}

Let us consider the complex flag manifold $G/B$, where $G$ is a
connected, simply connected,  simple, complex Lie group and
$B\subset G$ a Borel subgroup. Let $T$ be a maximal torus of a
compact real form of $G$, $\t$ its Lie algebra and $\Phi \subset
\t^*$ the corresponding set of roots. Consider an arbitrary
$W$-invariant inner product $\langle \ , \  \rangle$ on $\t$. The
Weyl group $W$ can be realized as the subgroup of the orthogonal
group of $(\t, \langle \ , \  \rangle)$ which is generated by the
reflections about the hyperplanes $\ker \alpha $,
 $\alpha\in \Phi^+$.
 To any root $\alpha$ corresponds the coroot
$$\alpha^{\vee}:= \frac{2\alpha}{\langle \alpha,
\alpha\rangle}$$ which  is an element of $\t$, by using the
identification of $\t$ and $\t^*$ induced by $\langle \ , \
\rangle$. If $\{\alpha_1, \ldots ,\alpha_l\}$ is a system of
simple roots then $\{\alpha_1^{\vee},\ldots, \alpha_l^{\vee}\}$ is
a  system of simple coroots.  Consider  $\{\lambda_1 ,\ldots ,
\lambda_l\} \subset \t^*$ the corresponding system of fundamental
weights, which are defined by
$\lambda_i(\alpha_j^{\vee})=\delta_{ij}$.

Let us recall the presentation of the cohomology\footnote{Only
cohomology with {\it real} coefficients will be considered
throughout this paper.} ring of $G/B$, as obtained by Borel in
\cite{B}.  First of all, one can assign to any weight $\lambda \in
\t^* $ a group homomorphism $T\to S^1$; the latter can be
 extended canonically to a group homomorphism $B \to \bC^*$ and
gives rise in this way to the complex line bundle  $L_{\lambda} =
G \times_B \bC $ over $G/B$.
 One shows that the ring homomorphism $S(\t^*)\to H^*(G/B)$
induced by $\lambda_i\mapsto c_1(L_{\lambda_i}),$ $1\leq i \leq
l$, is surjective; moreover it induces the  ring isomorphism
$$ H^*(G/B) \simeq \bR[\{\lambda_i\}]/I_W,$$ where $I_W$ is the
ideal of
$S(\t^*)=\bR[\lambda_1,\ldots,\lambda_l]=\bR[\{\lambda_i\}]$
generated by the $W$-invariant polynomials of strictly positive
degree.  We identify $H^*(G/B)$ with Borel's presentation and
denote them both by ${\mathcal H}$. So
$$\H=H^*(G/B)=\bR[\{\lambda_i\}]/I_W,$$ where $c_1(L_{\lambda_i})$ is identified with
the coset $[\lambda_i]$ of $\lambda_i$, $1\le i\le l$.  There are
two more things we would like to recall here:
 \begin{itemize}
 \item[-] by a result of Chevalley [C], there exist $l$
homogeneous, functionally independent polynomials $u_1,\ldots, u_l
\in S(\t^*)$, which generate $I_W$; \item[-] on $\H$ there exists
a natural inner product $( \ , \  )$, induced by the Poincar\'e
pairing.
\end{itemize}

Let us consider now the Hamiltonian system of Toda lattice type,
which consists of the standard symplectic manifold $(\bR^{2l},
\sum_{i=1}^l dr_i\wedge ds_i)$ with the Hamiltonian function
\begin{equation}\label{2}E(\{r_i\},\{s_i\})=\sum_{i,j=1}^l \langle \alpha_i^{\vee},
\alpha_j^{\vee}\rangle r_i r_j +\sum_{i=1}^l e^{2s_i}.
\end{equation} The following result, proved by Goodman and Wallach
in \cite{G-W}, gives details concerning the integrals of motion of
this system (note that the latter is completely integrable):

\begin{theorem}\label{gw1} {\rm (see \cite{G-W})} There exist $l$
functionally independent functions $$E=\tilde{F}_1, \tilde{F}_2,
\ldots, \tilde{F}_l :\bR^{2l}\to \bR$$ each of them uniquely
determined by:

\begin{itemize}
 \item[(i)] $\tilde{F}_k(\{r_i\},\{s_i\}) = F_k(\{e^{2s_i}\},\{r_i\})$, where
 $F_k$ is a polynomial in variables
 $e^{2s_1},\ldots,e^{2s_l}$, $r_1,\ldots,r_l,$ homogeneous with respect to $e^{s_1},\ldots ,
e^{s_l}, r_1,\ldots, r_l$;

 \item[(ii)] $\{\tilde{F}_k,E\}=0$, where $\{ \ , \ \}$ denotes the
 Poisson bracket of functions on $\bR^{2l}$;

 \item[(iii)] $F_k( 0, \ldots, 0, \lambda_1,\ldots, \lambda_l)=
u_k( \lambda_1,\ldots, \lambda_l)$ as elements of $S(\t^*)$.
\end{itemize}
\end{theorem}

Consider  now the formal multiplicative variables $q_1,\ldots,q_l$
which are assigned degree $4$ (note that the coset of $\lambda_j$
in $\bR[\{\lambda_i\}]/I_W$, which is the same as
$c_1(L_{\lambda_j})$ in $H^*(G/B)$, has degree 2). Occasionally,
$q_i$ will stand for $e^{t_i}$, $1\leq i \leq l$, where
$t_1,\ldots , t_l$ are real numbers, so that the differential
operators $\frac{\partial}{\partial t_i}$ on $\H \otimes
\bR[\{q_i\}]$ will be well defined.

Our goal is to  prove the following result:

\begin{theorem}\label{main} Let $\circ$ be an $\bR[\{q_i\}]$-linear
product on $\H\otimes \bR[\{q_i\}]$ with the following properties:
\begin{itemize}
\item[(i)] $\circ$ preserves the graduation induced by $\deg
[\lambda_i] =2$ and $\deg q_i=4$;
\item[(ii)] $\circ$ is a
deformation of the usual product, in the sense that if we formally
replace all $q_i$ by $0$, we obtain the usual product on $\H$;
\item[(iii)] $\circ$ is commutative;
\item[(iv)] $\circ$ is associative;
\item[(v)] $\sum_{i,j=1}^l\langle \alpha_i^{\vee}, \alpha_j^{\vee} \rangle
[\lambda_i]\circ [\lambda_j] =\sum_{i=1}^l \langle
\alpha_i^{\vee}, \alpha_i^{\vee} \rangle q_i$;
\item[(vi)] $\frac{\partial}{\partial
t_i}([\lambda_j] \circ a)=\frac{\partial}{\partial
t_j}([\lambda_i] \circ a)$, for
 any $a\in \H$, $1\leq i,j\leq l$.
\end{itemize}
 Then the ring $(\H \otimes
\bR[\{q_i\}], \circ)$ is generated by $[\lambda_1],\ldots,
[\lambda_l], q_1, \ldots, q_l$, subject to the relations
\begin{equation}\label{qrelation}F_k(\{-\langle \alpha_i^{\vee},\alpha_i^{\vee}\rangle
q_i\} ,\{[\lambda_i]\circ\})= 0, \end{equation} $1\leq k\leq l,$
where the polynomials $F_k$ are given by Theorem 1.1.
\end{theorem}

Our proof is purely algebraic, but the following geometric ideas
stay behind it: We assign to $\circ$ the 1-form $\omega$ on
$H^2(G/B)$ with values in ${\rm End}\H$, given by
$$\omega_t(X)(Y)= X\circ Y,$$ where $t=t_1[\lambda_1]+\ldots
+t_l[\lambda_l]\in H^2(G/B)$ and $X,Y\in \H$ (the convention
$q_i=e^{t_i}$ is in force). Consider the Dubrovin type connection
$\nabla^h = d +\frac{1}{h}\omega$ (cf. [D]) on the vector bundle
$\H \times H^2(G/B)\to H^2(G/B)$. Conditions (iv) and (vi) say
that $\nabla^{h}$ is a flat connection, for all $h\ne 0$. Let $( \
, \ )$ denote the Poincar\'e pairing on $\H$. We are able to
construct parallel sections $s:H^2(G/B)\to \H$ of the  connection
dual to $\nabla^h$, i.e. the one corresponding to $\omega^T$,
where
$$(\omega_t^T(X)(Y), Z)=(Y, \omega_t(X)(Z)),$$ $X,Y,Z\in \H$.
More precisely, we
find certain ``formal" solutions $s$ of the system
$$h\frac{\partial s}{\partial t_i} = \omega_t^T([\lambda_i])(s),$$
$1\le i\le l$ (for the details, see section 4). The main
difficulty is to show that the integrals of motion of the quantum
Toda lattice are {\it quantum differential operators} for $\circ$,
i.e. they vanish all functions $(s,1) : H^2(G/B)\to \bR$, where
$s$ is a parallel section as before: by results of Givental [G]
(see also [CK, section 10.3]), such differential operators induce
relations, and it is not difficult to see that those relations are
just (\ref{qrelation}). Now from condition (v) we can deduce that
the degree 2 integral of motion
--- call it $H$
--- {\it is} a quantum differential operator. Because $H$ commutes
with any other integral of motion, the latter is also a quantum
differential operator (this idea has also been used by B. Kim in
[K]).

{\bf Remarks.} 1. We only have to show that the relations
(\ref{qrelation}) hold in $(\H \otimes \bR[\{q_i\}], \circ)$: by a
general result of Siebert and Tian \cite{S-T}, they generate the
whole ideal of relations.

 2. Properties of the three-point Gromov-Witten invariants
 $\langle \ | \ | \ \rangle_d$
 (see for instance Fulton-Pandharipande [F-P]) show
that the hypotheses of Theorem \ref{main}  are satisfied by the
quantum product $\star$ on the (small) quantum cohomology ring of
$G/B$: Condition (v) follows immediately from the equation
\begin{equation}\label{lstar}[\lambda_i]\star [\lambda_j] = [\lambda_i][\lambda_j]
+\delta_{ij}q_j,\end{equation} and the fact that the degree two
homogeneous generator of $I_W$ (see the mention to Chevalley's
result from above) is
$$u_1 = \sum_{i,j=1}^l \langle
\alpha_i^{\vee},\alpha_j^{\vee}\rangle \lambda_i\lambda_j.$$ As
about (\ref{lstar}), it can be proved in an elementary way (see
\cite{K} or \cite{M1}). Condition (vi) is a direct consequence of
the definition
$$[\lambda_i]\star a = \sum_{d=(d_1,\ldots, d_l)\ge 0}
([\lambda_i] \star a)_d q_1^{d_1}\ldots q_l^{d_l} \  {\rm with} \
(([\lambda_i]\star a)_d, b) =\langle [\lambda_i] | a | b \rangle_d
\ {\rm for  \ all } \ b\in H^*(G/B)$$ and the ``divisor property"
$$\langle [\lambda_i] | a | b \rangle_d  = d_i \langle a |b \rangle_d.$$
 We recover in this way Kim's result
on $QH^*(G/B)$ (see \cite{K}). The main achievement of our paper
is that it shows that Kim's presentation of $QH^*(G/B)$ can be
deduced in an elementary way, by using very few of the properties
of the quantum product $\star$. For instance, the Frobenius
property
$$(a\star b, c)=(a, b\star c), \qquad a,b,c \in \H$$
is not needed in our proof.

3. In [M2] we constructed the ``combinatorial" quantum cohomology
ring and then we used Theorem \ref{main} in order to prove that
its isomorphism type is the  one expected by the theorem of Kim.

4. The main result of [M3] is an extension of Theorem \ref{main}:
we were able to obtain a similar connection between the small
quantum cohomology of the {\it infinite dimensional} generalized
flag manifold and the integrals of motion of the {\it periodic}
Toda lattice.

{\bf Acknowledgements.} I would like to thank Martin Guest and
Takashi Otofuji for several discussions on the topics of the
paper. I would also like to thank Roe Goodman for an illuminating
correspondence concerning Toda lattices. I am grateful to Lisa
Jeffrey for a careful reading of the manuscript and for suggesting
several improvements. I am also thankful to the referee for many
valuable suggestions.

\section{ Toda lattices according to Goodman and Wallach}
The goal of this section is to present two results of Goodman and
Wallach\cite{G-W}, which will be essential ingredients for the
proof of Theorem \ref{main}. Let us consider the $(ax+b)$-algebra
corresponding to the coroot system of $G$. By definition, this is
 the Lie algebra  $$(\b=\t^* \oplus \u, [ \ ,\ ]),$$ where $\u$
 has a basis $X_1,\ldots,X_l$ such that:
\begin{equation}\label{commutation}[\lambda_i,\lambda_j]=0, \quad
[\lambda_i,X_j]=\delta_{ij}X_j, \quad [X_i,X_j]=0,\end{equation}
$1\leq i,j\leq l$. The set $S(\b)$ of polynomial functions on
$\b^*$ is a Poisson algebra and by (\ref{commutation}) we have
$$\{\lambda_i,\lambda_j\}=0, \quad \{\lambda_i,X_j\}=\delta_{ij}X_j, \quad
\{X_i,X_j\}=0,$$ $1\leq i,j\leq l$.

On the other hand, one can easily see that the Poisson bracket of
functions on the standard symplectic manifold $(\bR^{2l},
\sum_{i=1}^l dr_i\wedge ds_i)$ satisfies
$$\{r_i,r_j\}=0, \quad \{r_i, e^{s_j}\}=\delta_{ij}e^{s_j}, \quad
\{e^{s_i}, e^{s_j}\}=0,$$ $1\leq i,j\leq l$. We deduce that the
Poisson subalgebra $\bR[e^{s_1},\ldots , e^{s_l}, r_1, \ldots,
r_l]$ of $C^{\infty}(\bR^{2l})$ is  isomorphic to $S(\b)$ via
\begin{equation}\label{transformation} X_i \mapsto e^{s_i}, \quad \lambda_i \mapsto r_i,
\end{equation} $1\leq i \leq l$. In this way, integrals of motion
of the Hamiltonian system determined by (\ref{2}) can be obtained
from elements of the space $S(\b)^{\{, \}}$, which is the $\{,
\}$-commutator in $S(\b)$ of the polynomial
\begin{equation}\label{hamiltonian} \sum_{i,j=1}^l \langle
\alpha_i^{\vee}, \alpha_j^{\vee}\rangle \lambda_i \lambda_j
+\sum_{i=1}^l X_i^2. \end{equation}

Let us consider now the universal enveloping algebra $$U(\b)=T
(\b)/\langle x\otimes y-y\otimes x -[x,y], x,y \in \b\rangle$$
with the canonical filtration $\{0\}=U_0(\b)\subset U_1(\b)
\subset \ldots $ (see e.g. [H, section 17.3]). We say that an
element $f$ of $U(\b)$ has {\it degree} $m$ if $m$ is the smallest
positive integer with the property that $f\in U_m(\b)$. There
exists a vector space isomorphism
$$\phi:S(\b) \to U(\b)$$
induced by the symmetrization map followed by the canonical
projection (see [H, Corollary E, section 17.3]). Since $\t^*$ and
$\u$ are abelian, the element of $S(\b)$ described by
(\ref{hamiltonian}) is mapped by $\phi$ to $$\Omega:=
\sum_{i,j=1}^l \langle \alpha_i^{\vee}, \alpha_j^{\vee}\rangle
\lambda_i \lambda_j +\sum_{i=1}^l X_i^2,$$ the right hand side
being regarded this time as an element of $U(\b)$.

 The complete integrability of the Toda lattice follows from the
 following two theorems of Goodman and Wallach:
\begin{theorem}\label{gw2} {\rm (see \cite{G-W})} The Poisson bracket commutator
$S(\b)^{\{,\}}$ is mapped by $\phi$ isomorphically onto the space
$U(\b)^{[,]}$ of all $f\in U(\b)$ with the property that
$[f,\Omega]=0$.
\end{theorem}

\begin{theorem}\label{gw3} {\rm (see \cite{G-W})} The map
 $\mu: U(\b)\to U(\t^*)=S(\t^*)$
induced by the natural Lie algebra homomorphism $\b\to \t^*$
establishes an algebra isomorphism between $U(\b)^{[,]}$ and the
ring $S(\t^*)^W$ of $W$-invariant polynomials. Hence there exist
$\Omega=\Omega_1, \Omega_2, \ldots, \Omega_l \in U(\b)$, each of
them uniquely determined by
\begin{itemize}
\item[(i)] $[\Omega_k, \Omega]=0,$
\item[(ii)] $\mu(\Omega_k)=u_k$ and $\deg\Omega_k = \deg u_k$.
\end{itemize}
Moreover, $\Omega_k$ is contained in the subring of $U(\b)$ which
is spanned by the elements of the form $ X^{2I}\lambda^J$.
\end{theorem}

\noindent {\bf Remark.} The integrals of motion of the Toda
lattice mentioned in
 Theorem \ref{gw1}
 are obtained from $\phi^{-1}(\Omega_k)$, $1\le k \le l$ by the transformations
(\ref{transformation}): by the last statement of Theorem
\ref{gw3}, the result are  polynomial expressions in variables
$e^{2s_1}, \ldots, e^{2s_l}, r_1,\ldots, r_l$ and these are what
we denoted $F_k(e^{2s_1}, \ldots, e^{2s_l}, r_1,\ldots, r_l)$.

\section{ Relations in $(\H\otimes \bR[q_i], \circ)$}

Consider the representation $\rho$ of $\b$ on $C^{\infty}(\bR^l)$
given by:
$$\rho (\lambda_i)=2 \frac{\partial}{\partial t_i}, \ \
\rho (X_i)=\frac{2\r1}{h} \sqrt{\langle \alpha_i^{\vee} ,
 \alpha_i^{\vee}\rangle}
e^{\frac{t_i}{2}}\cdot ,$$ $1\leq i \leq l$, where $h$ is a
nonzero real parameter. The differential operators
 $$D_k=h^{ \deg\Omega_k}\rho(\Omega_k),$$ $1\leq k \leq l$,
 will be the crucial objects of the proof of Theorem 1.2.

Since $F_k$ is homogeneous in variables $e^{s_i}, r_i$, it follows
that $\Omega_k$ --- being  obtained from $F_k$ after applying
$\phi$, up to the replacements (\ref{transformation})
--- has a presentation as a homogeneous, symmetric polynomial in
the variables $X_i, \lambda_i$. We use the commutation relations
(\ref{commutation}) in order to express $\Omega_k$ as a linear
combination of elements of the form $ X^{2I}\lambda^J$ (see
Theorem 2.2).
 The polynomial expression we
obtain in this way appears as
$$ \Omega_k=F_k(\{X_i^2\}, \{\lambda_i\}) +f_k(\{X_i^2\},\{ \lambda_i\})$$
where $$\deg f_k < \deg F_k.$$ Consequently  $D_k$ appears as a
polynomial expression  $D_k(e^{t_1},\ldots ,
e^{t_l},h\frac{\partial}{\partial t_1}, \ldots ,
h\frac{\partial}{\partial t_l}, h)$, the last ``variable", $h$,
being due to the possible occurrence  of $f_k$.

 Amongst all $D_k$, $1\leq k \leq l$, the operator $D_1=h^2\rho(\Omega)$
plays a privileged role, and we write
\begin{equation}\label{h}
H:=\frac{1}{4}D_1=h^2\sum_{i,j=1}^l\langle \alpha_i^{\vee},
\alpha_j^{\vee} \rangle \frac{\partial ^2}{\partial t_i \partial
t_j}- \sum_{j=1}^{l} \langle \alpha_j^{\vee}, \alpha_j^{\vee}
\rangle e^{t_j}. \end{equation}

Below we will see that the polynomial $D_k(\{Q_i\}, \{\Lambda_i\},
h)\in \bR[Q_1, \ldots, Q_l, \Lambda_1, \ldots, \Lambda_l, h]$
obtained from
 $D_k$ by the replacements $e^{t_i}\mapsto Q_i$, $h\frac{\partial}{\partial
t_i} \mapsto \Lambda_i$,
  $1\leq i \leq l$,  satisfies the hypotheses of the following theorem.

\begin{theorem}\label{criterion} Let $\circ$ be a product on $\H\otimes \bR[\{q_i\}]$
with the properties (i)-(vi) from Theorem \ref{main}. Suppose that
$D=D(\{Q_i\}, \{\Lambda_i\}, h)\in \bR[\{Q_i\}, \{\Lambda_i\},h]$
satisfies
\begin{itemize}
\item[(a)] $[D(\{e^{t_i}\}, \{h\frac{\partial}{\partial t_i}\},h),
H(\{e^{t_i}\},\{h\frac{\partial}{\partial t_i}\}, h)]=0$,
\item[(b)] $D(0,
\ldots, 0, \Lambda_1, \ldots, \Lambda_l, h)$ does not depend on
$h$,
\item[(c)] $D(0,\ldots , 0, \lambda _1, \ldots ,
\lambda_l,0)\in S(\t^*)^W$.
\end{itemize}
 Then the relation $D(\{q_i\}, \{[\lambda_i]\circ\},  0)=0$ holds in the ring
 $(\H\otimes \bR[\{q_i\}],\circ)$.
\end{theorem}

The proof of  this theorem will be done in the next section. Now
we will show  how can be used Theorem \ref{criterion} in order to
prove the main result of the paper.

\noindent {\it Proof of Theorem \ref{main}}.  By Remark 1 in the
introduction, we only have to show that the relations
(\ref{qrelation}) hold for all $1\le k \le l$. To this end we note
that $D=D_k$ satisfies the hypotheses of  Theorem \ref{criterion}:
 (a) follows from the fact that $\rho$ is  a Lie algebra
representation, and (b) and (c) from Theorem \ref{gw3} (ii). We
obtain the relation $D_k(\{q_i\}, \{[\lambda_i]\circ\},  0)=0$,
which is just (\ref{qrelation}). \qed

\section{Proof of Theorem \ref{criterion}}

 Let us begin by picking a
basis of $\H$ which consists of homogeneous elements (e.g. the
Schubert basis): this will allow us to identify $\H$ with $\bR^n$,
where $n=\dim \H$, and the endomorphism $[\lambda_i]\circ$ of $\H$
with an element $B_i$ of the space $M_n(\bR[e^{t_j}])$ of $n\times
n$ matrices whose coefficients are polynomials in $e^{t_1},
\ldots, e^{t_l}$. Let $\circ$ be a product which satisfies the
hypotheses of Theorem \ref{main}.

\begin{lemma}\label{lemma} Fix $i\in \{1,\ldots, l\}$ and take $a\in \H$. Write
\begin{equation}\label{da}[\lambda_i] \circ a=\sum_{d=(d_1, \ldots
, d_l)\geq 0} ([\lambda_i] \circ a)_dq^d
\end{equation} with $([\lambda_i] \circ a)_d \in \H$. If
$d=(d_1,\ldots,d_l)\ne 0$ such that $([\lambda_i]\circ a)_d\neq
0$, then $d_i\neq 0$. In other words, any non-zero term in the
right hand side of (\ref{da}) which is different from
$([\lambda_i] \circ a)_0=[\lambda_i]  a$ must be a multiple of
$q_i$.
\end{lemma}

\begin{proof} Condition (vi) from Theorem \ref{main} reads
$$\frac{\partial}{\partial t_i}B_j = \frac{\partial}{\partial t_j}B_i.$$
Hence there exists
 $M\in M_n(\bR[e^{t_j}])$  such that
$$B_i=B_i' + \frac{\partial}{\partial t_i} M,$$
where $B_i'$ is constant, for any $1\leq i \leq l$. It remains to
notice that the derivative with respect to $t_i$ of a monomial in
$e^{t_1}, \ldots, e^{t_l}$ contains only nonzero powers of
$e^{t_i}$, or else it is 0.
\end{proof}

 As pointed out in the introduction, $\H$ has
a natural inner product $( \ , \ )$, namely the Poincar\'e
pairing. Denote by $([\lambda_i]\circ)^T$ the endomorphism of $\H$
which is transposed to $[\lambda_i] \circ$ with respect to this
product, i.e.
$$([\lambda_i] \circ a, b)=(a,([\lambda_i]\circ)^Tb ), \quad a,b\in \H.$$
 Also denote by $A_i$ the  matrix of  $([\lambda_i]\circ)^T$ with respect to the basis
of $\H$ which is the dual with respect to $( \ , \ )$ of our
original basis: of course $A_i$ coincides  with the transposed of
the matrix of $[\lambda_i]\circ$ with respect to the original
basis. Now we  want the ordering of the original basis of $\H$ to
be decreasing with respect to the degrees of its elements. From
condition (i) from Theorem \ref{main} and Lemma \ref{lemma} it
follows that for any $i\in \{1,\ldots, l\}$, the matrix $A_i$ can
be decomposed as
$$A_i=A'_i+A''_i(e^{t_j})$$ where $A'_i$ is strictly lower
triangular and its coefficients do not depend on $t$ and $A_i''$
is strictly upper triangular,  its coefficients being linear
combinations of
$$e^{td}:=e^{t_1d_1} \ldots e^{t_l d_l},$$
 where $d_1,\ldots, d_l$ are nonnegative integer numbers  with $$d_i>0.$$

Consider the  PDE system:
\begin{equation}\label{pde}
h\frac{\partial}{\partial t_i}s =([\lambda_i]\circ)^T ( s),
\end{equation} $ 1\leq i \leq l,$ where the map  $s=s(t_1, \ldots ,t_l)$
takes values in $\H$ and $h$ is a nonzero real parameter. Some
algebraic formalism is needed in order to provide solutions to
(\ref{pde}). Let $\R$ be an arbitrary commutative, associative
real algebra with unit. For $V=\R^n$ or $V=M_n(\R)$ we denote by
$$V[t_i][[e^{t_i}]]:=V\otimes \R[t_1,\ldots ,t_l][[e^{t_1},\ldots ,e^{t_l}]]$$
the space of formal series  $$f=\sum_{d=(d_1, \ldots , d_l)\geq
0}f_d e^{td}$$ where $f_d$ is a polynomial in variables $t_1,
\ldots , t_l$ with coefficients in $V$. The operator
$\frac{\partial}{\partial t_i}$  acts in a natural way on
$V[t_i][[e^{t_i}]]$ via $$\frac{\partial}{\partial
t_i}(f_de^{td})= (\frac{\partial f_d}{\partial t_i}
+d_if_d)e^{td}.$$ We use the same formula
\begin{equation}\label{action}
(\sum_{d\ge 0} f_de^{td})(\sum_{d\ge 0} g_de^{td})=\sum_{d\ge
0}(\sum_{d_1+d_2=d}f_{d_1}g_{d_2})e^{td}\end{equation} in order to
define both:
\begin{itemize}
\item[-] an action of $M_n(\R)[t_i][[e^{t_i}]]$ on
$\R^n[t_i][[e^{t_i}]]$ (take $f_d\in M_n(\R)[t_i]$, $g_d\in
\R^n[t_i]$);
\item[-] a multiplication on $M_n(\R)[t_i][[e^{t_i}]]$ (take $f_d,
g_d\in M_n(\R)[t_i]$).
\end{itemize}  Alternatively, we can use the ring structure of
$\R[t_i][[e^{t_i}]]$ induced by the same formula (\ref{action})
(take $f_d,g_d \in \R[t_i]$), the identifications
$$M_n(\R)[t_i][[e^{t_i}]]=M_n(\R[t_i][[e^{t_i}]]), \ \R^n[t_i][[e^{t_i}]]=
(\R[t_i][[e^{t_i}]])^n$$ and the usual matrix multiplication
rules.

Our aim is to find solutions $s$ of the system (\ref{pde}) in the
space $\H[t_i][[e^{t_i}]]= \bR^n[t_i][[e^{t_i}]]$, where $\H$ has
been identified with $\bR^n$ via the basis  which is the dual with
respect to $( \ , \ )$ of the original basis (see above). The
following result will help us to this end:

\begin{proposition}\label{prop} Let   $\A_1, \ldots ,\A_l\in M_n(\R[e^{t_i}])$ be matrices which
satisfy:
\begin{itemize}
 \item[(a)]  $\A_i$ commutes with $\A_j$ for any two
 $i,j$;
\item[(b)] $\frac{\partial}{\partial t_i}\A_j=\frac{\partial}{\partial
t_j}\A_i$ for any two
 $i,j$;
\item[(c)] for any $i\in \{1,\ldots, l\}$ we can decompose $\A_i$ as
$$\A'_i+\A''_i(e^{t_j})$$ where $\A_i'$ is strictly lower
triangular and its coefficients do not depend on $t$, and $\A''_i$
is strictly upper triangular,
 its coefficients being linear combinations of
$e^{td}:=e^{t_1d_1} \ldots e^{t_l d_l}$, where $d_1 , \ldots ,
d_l$ are nonnegative integer numbers with
$$d_i>0.$$
\end{itemize}
 Consider the
 PDE system
\begin{equation}\label{system} \frac{\partial g}{\partial t_i} =\A_ig, \end{equation} $1\leq i
\leq l$, where
 $g=\sum_{d\geq 0}g_de^{td}
\in \R^n[t_i][[e^{t_i}]].$ The system has a unique solution $g$
with $g_0^0$ (the constant term of the polynomial $g_0\in
\R^n[t_i]$) prescribed.
\end{proposition}

 The following elementary lemma will be needed in the proof:

\begin{lemma}\label{lemma2}  Let  $\A\in M_n(\R)$ be a matrix and $g\in\R ^n[t]$ a
polynomial. Consider the differential equation:{\rm
$$\frac{\text{d}f}{{\text d}t}=\A f +g,$$}
where $f$ is in $\R ^n[t]$.

 (i) If $\A$ is invertible, then  we
have a unique solution $f$.

 (ii) If $\A$ is nilpotent, then  the equation has a unique solution $f$ with the
 constant term $f_0\in \R^n$ prescribed.
 \end{lemma}

\begin{proof} Put $g=\sum_{k=0}^p g_k t^k$ and look for $f$ as
$\sum_{j=0}^m f_j t^j$, where $g_k, f_j\in \R^n$. The  proof is
straightforward.
\end{proof}

\noindent {\it Proof  of Proposition \ref{prop}}. We will  prove
this result by induction on $l\geq 1$.  First take $l=1$ and solve
the equation
$$\frac{\text{d} g}{\text{d} t}=\A_1g,$$
where $g=\sum_{k\geq 0}g_{k}e^{tk}$, $g_{k}\in \R^n[t]$. Decompose
$\A_1$ as $\sum_{k\geq 0}(\A_1)^ke^{tk}$, where $(\A_1)^k\in
M_n(\R)$. Identify the coefficients of $e^{tk}$  and then
determine the polynomials $g_0, g_1, g_2 ,\ldots $ recursively by
Lemma \ref{lemma2}
 (notice that the matrix $(\A_1)^0=\A'_1$ is strictly lower triangular).

The induction step from $l-1$ to $l$ now follows. The idea is to
put $\S=\R[t_l][[e^{t_l}]]$ and note that we have
\begin{equation}\label{XX}\R^n[t_1,\ldots,t_l][[e^{t_1},\ldots,e^{t_l}]]=
\S^n[t_1,\ldots,t_{l-1}][[e^{t_1},\ldots,e^{t_{l-1}}]].\end{equation}
In other words, the $g$ we are looking for can be written as:
$$g=\sum_{d=(d_1, \ldots , d_{l})\geq 0}g_de^{t_1d_1 +\ldots + t_{l}d_{l}}=
\sum_{r=(r_1, \ldots ,r_{l-1})\geq 0}h_re^{t_1r_1+\ldots
+t_{l-1}r_{l-1}},$$ where $g_d\in\R^n[t_1,\ldots ,t_{l}]$ and
$h_r\in\R^n[t_l][[e^{t_l}]][t_1, \ldots
,t_{l-1}]=\S^n[t_1,\ldots,t_{l-1}]$. The identification given by
(\ref{XX}) maps $\A_i g$ to $\A_i h$, where the latter $\A_i$ is
regarded as an element of $M_n(\S[e^{t_1},\ldots,e^{t_{l-1}}])$,
$1\le i\le l$.

Our aim is to solve the system
$$\frac{\partial h}{\partial t_i}=\A_ih, \ 1\leq i \leq l-1,$$
where $h\in \S^n[t_1,\ldots ,t_{l-1}][[e^{t_1},\ldots
,e^{t_{l-1}}]]$. The elements $\A_1,\ldots ,\A_{l-1}$  of
$M_n(\S[e^{t_1},\ldots,e^{t_{l-1}}])$ satisfy  the conditions (a),
(b) and (c) (with $l-1$ instead of $l$). By the induction
hypothesis, we know that the solution of the latter PDE is
uniquely determined by the degree zero term $h_0^0$ of the
polynomial
 $h_0\in \S^n[t_1,\ldots ,t_{l-1}].$
We require that $h_0^0\in \S^n =\R^n[t_l][[e^{t_l}]]$ is the
solution of the equation
\begin{equation}\label{eq}\frac{\partial
h_0^0}{\partial t_l}=\A_l^0h_0^0 \end{equation} where $\A_l^0\in
M_n(\R[e^{t_l}])$ is the first term of the decomposition
 $\A_l=\sum_{r\geq 0}\A_l^re^{t_1r_1+\ldots + t_{l-1}r_{l-1}}$.

In order to be more precise, we write
$$h_0^0=\sum_{k\geq 0}f_k(t_l)e^{kt_l}$$
where $f_k(t_l)\in \R^n[t_l]$, $k\ge 0$, and then we identify the
coefficients of $e^{t_lk}$ in both sides of (\ref{eq}). One
obtains the following sequence of differential equations:
\begin{equation}\label{sequence}\frac{{\text d}f_k}{{\text
d}t_l}+kf_k=(\A_l^0h_0^0)_k= \sum_{u+v=k}(\A_l^0)^uf_v
\end{equation} where $(\A_l^0h_0^0)_k$ symbolizes the coefficient of
$e^{kt_l}$ in $\A_l^0h_0^0$ and $(\A_l^0)^u\in M_n(\R)$ is the
coefficient of $e^{t_lu}$ in $\A_l^0\in M_n(\R[e^{t_l}])$.

We solve the sequence (\ref{sequence}) of differential equations
by using Lemma \ref{lemma2}. First we write (\ref{sequence}) as:
$$\frac{{\text d}f_k}{{\text d}t_l}+kf_k=(\A_l^0)^0f_k+b,$$
where $b\in \R^n[t_l]$ depends only on $f_0,\ldots ,f_{k-1}$. The
matrix $(\A_l^0)^0\in M_n(\R)$ is obviously $\A_l'$ (see condition
(c)), hence it is strictly lower triangular. A simple recursive
procedure provide solutions: specifying  only $f_0^0=g_0^0$
determines first $f_0$ and then $f_1, f_2, \ldots $.

The only thing that remains to be  proved is that  the $g$  we
just constructed satisfies:
\begin{equation}\label{only}\frac{\partial g}{\partial t_l} =\A_lg.
 \end{equation} To this end, we notice first that
$$\frac{\partial }{\partial t_i}
 (\frac{\partial g}{\partial t_l} -\A_lg)=\A_i(\frac{\partial g}{\partial t_l} -\A_lg),$$ for all $1\leq i \leq l-1$.
Also $\frac{\partial g}{\partial t_l}-\A_lg$ can be written as
 $$\frac{\partial g}{\partial t_l}-\A_lg=
 \sum_{r=(r_1,\ldots ,r_{l-1})\geq 0}q_re^{t_1r_1+\ldots +t_{l-1}r_{l-1}},$$
 with $q_r\in \R^n[t_l][[e^{t_l}]][t_1, \ldots, t_{l-1}]$.
The degree zero term $q_0^0$ of $q_0$ is obviously $\frac{\partial
h_0^0}{\partial t_l}-\A_l^0h_0^0$.
 From the choice of $h_0^0$ it follows that $q_0^0=0$.
By the induction hypothesis,  $q_0^0$ determines $\frac{\partial
g}{\partial t_l}-\A_lg$ uniquely, hence the latter  is zero. \qed

 We apply Proposition \ref{prop} for $\R=\bR$ and $\A_i=\frac{1}{h}A_i$ and deduce:
\begin{corollary}\label{cor} For any $a\in \H$ there exists\footnote{
Proposition \ref{prop} also says that such an $s_a$ is unique, but
we do not need that.} $s_a\in
 \H[t_i][[e^{t_i}]]$ which is a solution of the system (\ref{pde})
 and satisfies the condition $(s_a)_0^0=a$.
 \end{corollary}

 There exists a $\bR[t_i]\me$-bilinear extension of the product $\circ$ to
$\H[t_i]\me$. Similarly, the intersection pairing $(\cdot, \cdot)$
can be extended to a $\bR[t_i]\me$-bilinear map
$$\H[t_i]\me \times \H[t_i]\me\rightarrow \bR[t_i]\me.$$ The differential
operator $D (\{e^{t_i}\},\{h \frac{\partial}{\partial t_i}\}, h)$
acts on $\bR[t_i]\me$ and $\H[t_i]\me$ in an obvious way. This
action plays an important role, as we can see in the following
lemma:

\begin{lemma}\label{lemma3} (i) Suppose that the differential operator
$D(\{e^{t_i}\},\{h\frac{\partial}{\partial t_i}\},
 h)$ satisfies
$$D.(s_a, 1)=0  \  {\it for \ all} \  a\in \H  {\it \ and \ all} \  h \neq 0,$$
where $s_a$ is given by Corollary \ref{cor}. Then we have the
relation $D(\{q_i\},\{ [\lambda_i]\circ\} ,0)=0 $.

 (ii) The following equation holds:
$$H.(s_a, 1)=0,$$ for all $a\in \H$, where the differential operator $H$ is given by
(\ref{h}).
\end{lemma}

\begin{proof} (i) For any $a\in \H$ and any $f\in \H\otimes \bR[e^{t_i}]$ we have that
$$\hbar\frac{\partial}{\partial t_i}(s_a, f)= (\hbar\frac{\partial}{\partial
t_i}s_a, f) + (s_a,\hbar\frac{\partial}{\partial t_i}f) =
(([\lambda_i]\circ)^T s_a, f) + (s_a,\hbar\frac{\partial}{\partial
t_i}f) = (s_a, ([\lambda_i]\circ+h\frac{\partial}{\partial
t_i})f).$$ We deduce that \begin{equation}\label{11}
D(e^{t_i},h\frac{\partial}{\partial t_i},
 h).(s_a,f)=
 (s_a, D(e^{t_1}, \ldots , e^{t_l},
[\lambda_1]\circ+h\frac{\partial}{\partial t_1}, \ldots
,[\lambda_l]\circ+h\frac{\partial}{\partial t_l},
 h).f).\end{equation}
Replacing $f$ by 1 and denoting $$\D=D(e^{t_1}, \ldots , e^{t_l},
[\lambda_1]\circ+h\frac{\partial}{\partial t_1}, \ldots
,[\lambda_l]\circ+h\frac{\partial}{\partial t_l},
 h).1,$$ we obtain
\begin{equation}\label{d}(\D, s_a)=0\end{equation} for all $a\in \H$.

For the rest of the proof, ``degree" will refer to the variables
$e^{t_1}, \ldots, e^{t_l}$. Note that $\D$ is an element of
$\H\otimes \bR[e^{t_i}]$. Decompose it as
$$\D = {\mathcal D}_0 + \D_1 + \ldots + \D_m,$$
where $\D_k\in \H\otimes \bR[e^{t_i}]$ denotes the sum of all
monomials of degree $k$, $0\le k \le m$. Recall that the degree
zero term of $s_a$ is the polynomial $(s_a)_0\in \H\otimes
\bR[t_1,\ldots, t_l]$, with $(s_a)_0^0=a$. The degree zero term of
$(\D, s_a)$ is $(\D_0, (s_a)_0)$. From the vanishing of the latter
we obtain that $$(\D_0, (s_a)_0^0) = (\D_0, a) =0,$$ for all $a\in
\H$, hence $\D_0=0$.

Also the sum of the terms of degree 1 in $(\D, s_a)$ is zero.
Since $\D_0=0$, this implies that $$(\D_1, (s_a)_0) = 0.$$ As
before, we deduce that $\D_1=0$. We continue this process and show
inductively that $\D_k =0$, for all $0\le k \le m$, hence
$$\D=0.$$ Now we let $h$ approach zero and deduce the desired
relation: $$D(q_i, [\lambda_i]\circ ,0)=0 .$$

(ii) When computing $H.(s_a, 1)$ we only need the fact that
$$\hbar ^2\frac{\partial ^2}{\partial t_i \partial t_j}(s_a, 1)=
(s_a, [\lambda _i]\circ [\lambda _j]), $$ which can be deduced
immediately from (\ref{11}). This implies that
$$H.(s_a,1) = (s_a,\sum_{i,j=1}^l\langle \alpha_i^{\vee}, \alpha_j^{\vee} \rangle
[\lambda_i]\circ [\lambda_j] -\sum_{i=1}^l \langle
\alpha_i^{\vee}, \alpha_i^{\vee} \rangle e^{t_i}) =0,$$ where we
have used condition (v) from Theorem \ref{main}.
\end{proof}

Another important step will be made by the following lemma:

\begin{lemma}\label{lemma4} {\rm (Kim's lemma, see \cite{K})}
Let\footnote{Here $d>0$ means $d\ge 0$ and $d\ne 0$.}
 $g=g_0 +\sum_{d > 0} g_d e^{td}\in
\bR[t_i]\me$ be a formal series with the properties $g_0=0$ and
$H.g=0.$ Then $g=0.$
\end{lemma}

\begin{proof} Suppose $g\neq 0$. Fix $d\in \bZ^l$, $d_i\ge 0$, with $g_d\neq 0$ and
$|d|:=\sum_{i=1}^l d_i>0$ minimal. From $H.g=0$ it follows that
\begin{equation}\label{zero}
 \sum_{i,j=1}^l\langle \alpha_i^{\vee}, \alpha_j^{\vee}
\rangle\frac{\partial ^2}{\partial t_i\partial t_j}(g_de^{td})=0.
\end{equation}  On the other hand, we have
$$\sum_{i,j=1}^l\langle \alpha_i^{\vee},
\alpha_j^{\vee}\rangle\frac{\partial ^2}{\partial t_i\partial t_j} (e^{td})=\sum_{i,j=1}^l\langle \alpha_i^{\vee},\alpha_j^{\vee}\rangle d_i d_j e^{td}=
 || \sum_{j=1}^l d_j \alpha _j^{\vee}||^2 e^{td}>0.$$
 Hence (\ref{zero}) is impossible.
  \end{proof}

And now we are in a position to prove Theorem \ref{criterion}:

\noindent {\it Proof of Theorem \ref{criterion}}. By Lemma
\ref{lemma3}, it is sufficient to show that
$$g:=D.(s_a, 1)=
\sum_{d\geq 0}g_de^{td}$$ equals zero. Taking into account
(\ref{11}), we have that
\begin{align} g&=(s_a, D(e^{t_1}, \ldots , e^{t_l},
[\lambda_1]\circ+h\frac{\partial}{\partial t_1}, \ldots
,[\lambda_l]\circ+h\frac{\partial}{\partial t_l},
 h).1)\nonumber
\\{}&=(s_a, D(0, \ldots , 0,[\lambda_1]\circ,\ldots ,  [\lambda_l]\circ,
h )+R), \nonumber\end{align}
 where $R \equiv 0$ mod $\{e^{t_i}\}$. Hence the polynomial
$g_0$ must be the same as  $$(s_a, D( 0, \ldots, 0,[\lambda_1],
\ldots , [\lambda_l], h))_0$$ (in our notation, the subscript $0$
indicates the constant term with respect
 to $\{e^{t_i}\}$). By conditions (b) and (c), the latter expression is zero, hence $g_0=0$.
It remains to notice that $H.g=0$ (which follows from
 $H. (s_a, 1)=0$
and $[D, H]=0$) and  apply
  Lemma \ref{lemma4}. \qed

\end{document}